\documentstyle[amscd,amssymb,verbatim,12pt]{amsart}

\theoremstyle{plain}
\newtheorem{Prop}{Proposition}[section]
\newtheorem{Thm}[Prop]{Theorem}
\newtheorem{Cor}[Prop]{Corollary}

\newtheorem{Lem}[Prop]{Lemma}

\theoremstyle{definition}
\newtheorem{Def}[Prop]{Definition}
\newtheorem{Question}[Prop]{\bf Question}
\theoremstyle{remark}
\newtheorem{Rem}[Prop]{Remark}
\newtheorem{Problem}[Prop]{\bf Problem}

\def\dim{\mathop{\roman{dim}}}
\def\int{\mathop{\roman{int}}}

\def\1{^{-1}}

\def\dim{\text{dim}}
\def\diam{\text{diam}}

\def\asdim{\text{asdim}}
\def\Lasdim {\text{l-asdim}}
\def\ANasdim {\text{asdim}_{AN}}
\def\rdim{\text{r-dim}}

\def\UU{{\mathcal U}}
\def\VV{{\mathcal V}}

\def\RR{{\mathbb R}}
\def\ZZ{{\mathbb Z}}
\def\dokaz{{\bf Proof. }}
\numberwithin{equation}{section}
\begin{document}
\title[
Hurewicz Theorem for Assouad-Nagata dimension]%
   {Hurewicz Theorem for Assouad-Nagata dimension
}

\author{N.~Brodskiy}
\address{University of Tennessee, Knoxville, TN 37996, USA}
\email{brodskiy@@math.utk.edu}

\author{J.~Dydak}
\address{University of Tennessee, Knoxville, TN 37996, USA}
\email{dydak@@math.utk.edu}

\author{M.~Levin}
\address{Departament of Mathematics, Ben-Gurion University of the
Nagev, P.O.B. 653, Beer-Sheba 84105, Israel. }
\email{mlevine@@math.bgu.ac.il}

\author{A.~Mitra}
\address{University of Tennessee, Knoxville, TN 37996, USA}
\email{ajmitra@@math.utk.edu}

\date{ May 17, 2006}


\thanks{ The second-named and third-named authors were partially supported
by Grant No.2004047  from the United States-Israel Binational Science
Foundation (BSF),  Jerusalem, Israel.
}

\begin{abstract} 
Given a function $f\colon X\to Y$ of metric spaces, its
{\it asymptotic dimension} $\asdim(f)$ is the supremum of $\asdim(A)$
such that $A\subset X$ and $\asdim(f(A))=0$.
Our main result is
\begin{Thm} \label{ThmAInAbstract}
$\asdim(X)\leq \asdim(f)+\asdim(Y)$ for any large scale uniform function $f\colon X\to Y$.
\end{Thm}

\ref{ThmAInAbstract}  generalizes a
result of Bell and Dranishnikov \cite{BellDranish A Hurewicz Type} in which
$f$ is Lipschitz and $X$ is geodesic.
We provide analogs of \ref{ThmAInAbstract}
for Assouad-Nagata dimension $\dim_{AN}$ and asymptotic Assouad-Nagata dimension
$\ANasdim$. In case of linearly controlled asymptotic dimension
$\Lasdim$ we provide counterexamples to three questions of Dranishnikov \cite{DranPrList}.

As an application of analogs of \ref{ThmAInAbstract} we prove
\begin{Thm}  \label{ThmBInAbstract}
If $1\to K\to G\to H\to 1$ is an exact sequence of groups and $G$ is finitely generated, then
$$\ANasdim (G,d_G)\leq \ANasdim  (K,d_G|K)+\ANasdim (H,d_H)$$
for any word metrics metrics $d_G$ on $G$ and $d_H$ on $H$.
\end{Thm}
  \ref{ThmBInAbstract} extends a
result of Bell and Dranishnikov \cite{BellDranish A Hurewicz Type} for asymptotic dimension.
\end{abstract}

\maketitle

\medskip
\medskip
\tableofcontents

\section{Introduction}\label{section Introduction}

The well-known Hurewicz Theorem for maps (also known as Dimension-Lowering Theorem, see \cite{Engel}, Theorem 1.12.4 on p.109) says $\dim(X)\leq \dim(f)+\dim(Y)$ if $f\colon X\to Y$ is a closed map
of separable metric spaces and $\dim(f)$ is defined as the supremum
of $\dim(f^{-1}(y))$, $y\in Y$. Bell and Dranishnikov \cite{BellDranish A Hurewicz Type}
proved a variant of Hurewicz Theorem for asymptotic dimension without defining
the asymptotic dimension of a function.
However, Theorem 1 of \cite{BellDranish A Hurewicz Type} may be restated
as $\asdim(X)\leq \asdim(f)+\asdim(Y)$, where $\asdim(f)$ is the smallest integer $n$
such that $\asdim(f^{-1}(B_R(y)))\leq n$ uniformly for all $R > 0$.
As an application it is shown in \cite{BellDranish A Hurewicz Type}
that $\asdim(G)\leq \asdim(K)+\asdim(H)$ for any exact sequence
$1\to K\to G\to H\to 1$ of finitely generated groups.
That inequality was extended subsequently by Dranishnikov and Smith \cite{DS}
to all countable groups.

\par The purpose of this paper is to generalize Hurewicz Theorem to variants of asymptotic dimension:
asymptotic Assouad-Nagata dimension and Assouad-Nagata dimension.
In the process we produce a much simpler proof than that in \cite{BellDranish A Hurewicz Type}
and a stronger result: the function is only assumed to be large scale uniform instead of Lipschitz,
and the domain is not required to be geodesic.

One of the main tools is Kolmogorov's idea used in his solution to Hilbert's 13th Problem. In dimension theory it is known as Ostrand Theorem. Another tool is reformulating
Gromov's \cite{Gro asym invar} definition of asymptotic dimension in terms of
$r$-components of spaces. That leads to a definition of
asymptotic dimension of a function in terms of double-parameter components, a concept
well-suited for Kolmogorov Trick.

\section{Ostrand theorem for asymptotic dimension}\label{section Ostrand for asymptotic dimension}

The aim of this section is to prove a variant of Ostrand's Theorem
for large scale dimensions. As an application we present a simple proof of
the Logarithmic Law for large scale dimensions.

Intuitively, a metric space is of dimension $0$ at scale $r$ if it can be represented as
a collection of $r$-disjoint and uniformly bounded subsets.
The following definition of Gromov \cite{Gro asym invar} defines
asymptotic dimension at most $n$ in terms of being represented, for each $r > 0$, as the union
of $n+1$ sets of dimension $0$ at scale $r$. That corresponds to the well-known
property of topologically $n$-dimensional spaces.

\begin{Def} \label{DefOfAsDimOfASpace}
A metric space $X$ is said to be of {\it asymptotic dimension}
at most $n$ (notation: $\asdim(X)\leq n$) if there is a function
$D_X\colon R_+\to R_+$
such that for all $r > 0$ there is a cover
$\UU=\bigcup\limits_{i=1}^{n+1}\UU_i$ of $X$ so that
each $\UU_i$ is $r$-disjoint (that means $d(a,b)\ge r$ for any two points
$a$ and $b$ belonging to different elements of $\UU_i$) and the diameter of elements of
$\UU$ is bounded by $D_X(r)$.
\par We refer to the function $D_X$ as {\it an $n$-dimensional control function} for $X$.
When discussing variants of asymptotic dimension it is convenient to
allow $D_X$ to assume infinity as its value at some range of $r$.
\end{Def}

\begin{Def} \label{DefOfAssouadLinearlyControlled}

A metric space $X$ is said to be of {\it Assouad-Nagata dimension}
(see \cite{Lang-Sch Nagata dim} and \cite{Brod-Dydak-Higes-MitraNA})
at most $n$ (notation: $\dim_{AN}(X)\leq n$) if it has
an $n$-dimensional control function $D_X$ that is a dilation
($D_X(r)=c\cdot r$ for some $c > 0$).

A metric space $X$ is said to be of {\it asymptotic Assouad-Nagata dimension}
at most $n$ (notation: $\ANasdim (X)\leq n$) if it has
an $n$-dimensional control function $D_X$ that is linear
($D_X(r)=c\cdot r+b$ for some $b, c \ge 0$).

A metric space $X$ is said to be of {\it linearly controlled asymptotic dimension}
at most $n$ (Dranishnikov \cite{DranPrList}, notation: $\Lasdim (X)\leq n$) if it has
an $n$-dimensional control function $D_X$ satisfying $D_X(r)=c\cdot r$ for some $c > 0$
and for all $r$ belonging to some unbounded subset of $R_+$.
Strictly speaking, the original definition of Dranishnikov \cite{DranPrList} is formulated in terms of
Lebesque numbers. However, just as in the case of asymptotic dimension,
it is equivalent to our definition.

A metric space $X$ is said to be of {\it microscopic Assouad-Nagata dimension}
at most $n$ if it has
an $n$-dimensional control function $D_X\colon R_+\to R_+\cup\infty$ that is a dilation near $0$:
$D_X(r)=c\cdot r$ for some $c > 0$ and all $r$ smaller than some positive number $M$,
$D_X(r)=\infty$ for all $r\ge M$.
\end{Def}

\begin{Def} \label{DefOfDimControlFunctionForLargerk}
Given a metric space $X$ and $k\ge n+1\ge 1$ an {\it $(n,k)$-dimensional control function} for $X$
is a function $D_X\colon R_+\to R_+$ such that for any $r > 0$ there is
a family $\UU=\bigcup\limits_{i=1}^{k}\UU_i$ satisfying the following conditions:

\begin{enumerate}
\item each $\UU_i$ is $r$-disjoint,
\item each $\UU_i$ is $D_X(r)$-bounded,
\item each element $x\in X$ belongs to at least $k-n$ elements of $\UU$ (equivalently, 
$\bigcup\limits_{i\in T}\UU_i$ is a cover of $X$ for every subset $T$ of $\{1,\ldots,k\}$
consisting of $n+1$ elements).

\end{enumerate}
 \end{Def}

The following result is an adaptation of a theorem by Ostrand
\cite{OstrandDimofMetriSpacesHilbert}.

\begin{Thm} \label{OstrandTypeThm}
If $D^{(n+1)}_X$ is an $n$-dimensional control function of $X$
and one defines a sequence of functions $\{D^{(i)}_X\}_{i\ge n+1}$ inductively by
$D^{(i+1)}_X(r)=D^{(i)}_X(3r)+2r$ for all $i\ge n+1$,
then each $D^{(k)}_X$ is an $(n,k)$-dimensional control
function of $X$ for all $k\ge n+1$.
\end{Thm}

\dokaz The proof is  by induction on $k$. The
case of $k=n+1$ is obvious. Suppose the result holds for some $k \ge n+1$.
Let $\UU=\bigcup\limits_{i=1}^{k}\UU_i$ be a family such that each $\UU_i$ is $3r$-disjoint,
each $\UU_i$ is $D^{(k)}_X(3r)$-bounded,
and each element $x\in X$ belongs to at least $k-n$ elements of $\UU$.
Define 
$\UU'_i$ to be the $r$-neighborhoods of elements of $\UU_i$ for $i\leq k$.
Notice elements of $\UU'_i$ are $(D^{(k)}_X(3r)+2r)$-bounded
and are $r$-disjoint.
Define
$\UU'_{k+1}$ as the collection of all sets of the form $\bigcap\limits_{s\in S}A_{s}\setminus \bigcup\limits_{i\notin S}\UU'_i$, where $S$ is a subset of $\{1,\ldots,k\}$ consisting of exactly
$k-n$ elements and $A_s\in \UU_s$. 

Notice that any element of $\UU'_{k+1}$ is contained in a single
element of some $\UU_j$. Thus
elements of each $\UU'_i$ are $(D^{(k)}_f(3r)+2r)$-bounded. 

Given two different sets $A=\bigcap\limits_{s\in S}A_{s}\setminus \bigcup\limits_{i\notin S}\UU'_i$, where $S$ is a subset of $\{1,\ldots,k\}$ consisting of exactly
$k-n$ elements and $A_s\in \UU_s$,
and $B=\bigcap\limits_{t\in T}B_{t}\setminus \bigcup\limits_{i\notin T}\UU'_i$, where $T$ is a subset of $\{1,\ldots,k\}$ consisting of exactly
$k-n$ elements and $B_t\in \UU_t$, we need to show $A$ and $B$ are $r$-disjoint.
It is clearly so if $S=T$, so assume $S\ne T$. If $a\in A$, $b\in B$, and $d(a,b)< r$,
then there is $s\in S\setminus T$ such that $a\in A_s$ prompting $b\in U\in \UU'_s$,
a contradiction.

Suppose $x\in X$ belongs exactly
to $k-n$ sets $\bigcup \UU'_i$, $i\leq k$, and let $S=\{i\leq k \mid x\in \bigcup \UU'_i\}$.
If $x\notin \bigcup \UU'_{k+1}$, then $x$ must belong to $\bigcup\UU'_j$ for some $j\notin S$,
a contradiction.
Thus each $x\in X$ belongs to at least $k+1-n$ elements of $\{\bigcup \UU'_i\}_{i=1}^{k+1}$.
\hfill $\blacksquare$

The product theorem for asymptotic dimension was proved in
\cite{DJ} using maps to polyhedra. 
The product theorem for Nagata dimension was proved in
\cite{Lang-Sch Nagata dim} using Lipschitz maps to polyhedra. Below
we use Theorem \ref{OstrandTypeThm} to give a simplified proof for all dimension theories.
Note the metric on $X\times Y$ is the sum of corresponding metrics on $X$ and $Y$.

\begin{Thm} \label{ProductThm}
If $X$ and $Y$ are metric spaces, then $$D(X \times Y) \le
D(X) +D(Y),$$
where $D$ stands for any of the following dimension theories: asymptotic dimension,
asymptotic  Assouad-Nagata dimension, Assouad-Nagata dimension, or microscopic Assouad-Nagata dimension.
\end{Thm}

\dokaz Let $D(X)=m$, $D(Y)=n$ and
let $k=m+n+1$. Pick $(m,k)$-dimension control function
$D_X$ of $X$ and $(n,k)$-dimensional control function $D_Y$ of $Y$, both of the correct type
(arbitrary, linear, dilation, or a dilation near $0$).
There are families $\{{\mathcal U}_i\}_{i=1}^k$ in $X$
and $\{{\mathcal V}_i\}_{i=1}^k$ in $Y$ that are $r$-disjoint and
bounded by $D_X(r)$ and $D_Y(r)$ respectively, that cover $X$ and $Y$ at
least $k-m$ times and $k-n$ times respectively. Then the family
$\{{\mathcal U}_i \times {\mathcal V}_i\}_{i=1}^k$ covers $X \times Y$,
as for any point $(x,y)$, $x$ is contained in sets from at least
$k-m=n+1$ families from $\{{\mathcal U}_i\}_{i=1}^k$ and $y$ is contained in sets from at least
$k-n=m+1$ families from $\{{\mathcal V}_i\}_{i=1}^k$, so there is at least one index
$j$ such that $x$ is covered by $\UU_j$ and $y$ is covered by $\VV_j$. The family
$\{{\mathcal U}_i \times {\mathcal V}_i\}_{i=1}^k$ is $r$-disjoint and
is bounded by $D_X(r)+D_Y(r)$.
\hfill $\blacksquare$

\section{Components, dimension, and coarseness}

In this section we replace the language of $r$-disjoint families
by the language of $r$-components. This language has the advantage of being portable
to functions. It also allows for simple proofs of known results \ref{GroupAndFGenSubgroups}-\ref{UnionThmForAsdim}.

\begin{Def} \label{DefOfrComponents}
Let $f\colon X\to Y$ be a function of metric spaces, $A$ is a subset
of $X$, and $r_X$, $r_Y$ are two positive numbers. 

$A$
is {\it $(r_X,r_Y)$-bounded} if for any points $x,x'\in A$ we have
$$d_X(x,x')\le r_X\quad \text{and}\quad
d_Y(f(x),f(x'))\le r_Y. $$

An {\it $(r_X,r_Y)$-chain} in $A$ is a sequence of points $x_1,\ldots,x_k$ in $A$
such that for every $i<k$ the set $\{x_i,x_{i+1}\}$ is $(r_X,r_Y)$-bounded.

 $A$ is {\it $(r_X,r_Y)$-connected} if for any points $x,x'\in
A$ can be connected in $A$ by an $(r_X,r_Y)$-chain.
\end{Def}

Notice that any subset $A$ of $X$ is a union of its {\it
$(r_X,r_Y)$-components} (the maximal $(r_X,r_Y)$-connected subsets
of $A$).

\begin{Def} \label{DefOfCoarseFunction}
Let $f\colon X\to Y$ be a function of metric spaces. $f$ is called {\it large scale uniform}
if there is function $c_f\colon R_+\to R_+$ such that 
$d_X(x,y) \leq r$ implies $d_Y(f(x),f(y)) \leq c_f(r)$.
The function $c_f$
will be called {\it a coarseness control function} of $f$.
\end{Def}

Notice $f$ is Lipschitz if and only if it has a coarseness control function that is a dilation.
$f$ is {\it asymptotically Lipschitz} if and only if it has a coarseness control function that is linear.

The following Lemma describes a useful case in which double parameter components coincide with single parameter components. The proof of this Lemma is an easy exercise.

\begin{Lem}\label{CoarsnessLipViaComponents}
 Let $f\colon X\to Y$ be a function of metric spaces.
$c_f\colon R_+\to R_+$ is a coarseness control function of $f$ if and only if for any subset $A$ of $X$
its $(r,c_f(r))$-components coincide with its $r$-components.
\end{Lem}

The following Lemma describes the way we construct a subset of
$X$ with all $(r_X,r_Y)$-components being
$(R_X,R_Y)$-bounded. The proof of this Lemma is an easy exercise.

\begin{Lem}\label{Lemma How we get double control}
Let $f\colon X\to Y$ be a function of metric spaces, $B$ be a subset of
$X$, and $A$ be a subset of $Y$. If all $r_X$-components of $B$
are $R_X$-bounded and all $r_Y$-components of $A$ are
$R_Y$-bounded then all $(r_X,r_Y)$-components of the set $B\cap
f^{-1}(A)$ are $(R_X,R_Y)$-bounded.
\end{Lem}

\begin{Lem} \label{ComponentsOfUnions}
Let $f\colon X\to Y$ be a function of metric spaces and $A,B$ be
subsets of $X$. Suppose that all $(r^A_X,r^A_Y)$-components of $A$
are $(R^A_X,R^A_Y)$-bounded and all $(r^B_X,r^B_Y)$-components of
$B$ are $(R^B_X,R^B_Y)$-bounded. If $R^B_X+2r^B_X<r^A_X$ and
$R^B_Y+2r^B_Y<r^A_Y$ then all $(r^B_X,r^B_Y)$-components of $A\cup
B$ are $(R^A_X+2r^A_X,R^A_Y+2r^A_Y)$-bounded.
\end{Lem}

\begin{pf}
Let $x=x_1,x_2,\dots,x_n=x'$ form an $(r^B_X,r^B_Y)$-chain in
$A\cup B$.

Notice that, if for some indices $i<j$ we have $x_k\in B$ for all
$i< k< j$, then $x_{i+1}$ and $x_{j-1}$ are in one
$(r^B_X,r^B_Y)$-component of $B$ and therefore
$d_X(x_{i+1},x_{j-1})\le R^B_X$ and $d_Y(f(x_{i+1}),f(x_{j-1}))\le
R^B_Y$.

If $x_i, x_j\in A$ and $x_k\in B$ for all $i< k< j$, then
$$d_X(x_{i},x_{j})\le
d_X(x_{i},x_{i+1})+d_X(x_{i+1},x_{j-1})+d_X(x_{j-1},x_{j})\le
r^B_X+R^B_X+r^B_X<r^A_X$$ and, similarly, $d_Y(f(x_{i}),f(x_{j}))\le
r^B_Y+R^B_Y+r^B_Y<r^A_Y$. Thus the points $x_i, x_j$ belong to the
same $(r^A_X,r^A_Y)$-component of $A$. This implies that all
points in the chain $x=x_1,x_2,\dots,x_n=x'$ belonging to $A$ are
in one $(r^A_X,r^A_Y)$-component of $A$.

Now let $x_s$ be the first point in the chain belonging to $A$ and
$x_t$ be the last point in the chain belonging to $A$. Then
$$d_X(x,x')\le$$
$$d_X(x_1,x_{s-1})+d_X(x_{s-1},x_s)+d_X(x_s,x_t)+d_X(x_t,x_{t+1})+d_X(x_{t+1},x_
{n})\le
$$ 
$$ \le R^B_X+r^B_X+R^A_X+r^B_X+R^B_X<R^A_X+2r^A_X.$$
Similarly, $d_Y(f(x),f(x'))\le
R^B_Y+r^B_Y+R^A_Y+r^B_Y+R^B_Y<R^A_Y+2r^A_Y.$
\end{pf}

\begin{Cor} \label{rComponentsOfBigUnions}
Let $f\colon X\to Y$ be a function of metric spaces and
$\{B_i\}_{i=1}^n$ be subsets of $X$. Suppose that for every $i$
all $(r^{(i)}_X,r^{(i)}_Y)$-components of $B_i$ are
$(R^{(i)}_X,R^{(i)}_Y)$-bounded. If, for every $i<n$,
$R^{(i+1)}_X+2r^{(i+1)}_X<r^{(i)}_X$ and
$R^{(i+1)}_Y+2r^{(i+1)}_Y<r^{(i)}_Y$, then all
$(r^{(n)}_X,r^{(n)}_Y)$-components of $\bigcup\limits_{i=1}^n B_i$ are
$(R^{(1)}_X+2r^{(1)}_X,R^{(1)}_Y+2r^{(1)}_Y)$-bounded.
\end{Cor}

\begin{pf}
By induction using Lemma~\ref{ComponentsOfUnions}.
\end{pf}

\begin{Prop} \label{DimControlOfAneighborhood}
Suppose $A$ is a subset of a metric space $X$, $m\ge 0$, and $R > 0$.
If $D_A$ is an $m$-dimensional control function of $A$, then
$D_B(x):=D_A(x+2R)+2R$ is an $m$-dimensional control function
of the $R$-neighborhood $B=B(A,R)$ of $A$.
\end{Prop}
\dokaz Given $r > 0$ express $A$ as $\bigcup\limits_{i=1}^{m+1}A_i$
such that $(r+2R)$-components of $A_i$ are $D_A(r+2R)$-bounded.
Given an $r$-component of $B_i:=B(A_i,R)$, each point
in that component is $R$-close to a single $(r+2R)$-component of $A_i$.
Therefore $r$-components of $B_i$ are $(D_A(r+2R)+2R)$-bounded.
\hfill $\blacksquare$

\begin{Def} \label{DefOfrDim}
Given a metric space $X$ and $r > 0$ the {\it $r$-scale dimension} $\rdim(X)$ is the smallest integer $n\ge 0$
such that $X$ can be expressed as $X=\bigcup\limits_{i=1}^{n+1}X_i$
and $r$-components of each $X_i$ are uniformly bounded.
\end{Def}
Notice $\asdim(X)$ is the smallest integer $n$ such that $\rdim(X)\leq n$ for all $r > 0$.
Also, $\asdim(X)$ is the smallest integer such that for each $r > 0$
the space $X$ can be expressed as $X=\bigcup\limits_{i=1}^{n+1}X_i$
so that $\rdim(X_i)\leq 0$ for each $i\leq n+1$.

\begin{Cor} \label{AsdimAndrDimUnion}
Suppose $X$ is a metric space. If, for every $r > 0$, there is a subspace $X_r$ of $X$
such that $\asdim(X_r)\leq n$ and $\rdim(X\setminus X_r)\leq n$,
then $\asdim(X)\leq n$.
\end{Cor}
\dokaz 
Express $X\setminus X_r$ as $\bigcup\limits_{i=1}^{n+1}A_i$
such that $r$-components of $A_i$ are $R$-bounded.
Express $X_r$ as $\bigcup\limits_{i=1}^{n+1}B_i$
so that $(R+2r)$-components of $B_i$ are $M$-bounded.
By Lemma~\ref{ComponentsOfUnions}, $r$-components of $A_i\cup B_i$ are $(M+2R+4r)$-bounded.
\hfill $\blacksquare$

\begin{Cor} \label{GroupAndSubgroups}
Suppose $G$ is a group with a left-invariant metric
and $G_r$ is the subgroup of $G$ generated by $B(1_G,r)$. If, for every $r > 0$, $\asdim(G_r)\leq n$,
then $\asdim(G)\leq n$.
\end{Cor}
\dokaz 
Consider $g_s\in G$, $s\in S$, so that $[g_s]$ enumerate $G/G_r\setminus [G_r]$.
Notice $Y=\bigcup\limits_{s\in S} g_s\cdot G_r$ is of $\rdim(Y)\leq \asdim(G_r)$
and $Y\cup G_r=G$.
\hfill $\blacksquare$

\begin{Cor}[Dranishnikov-Smith \cite{DS}] \label{GroupAndFGenSubgroups}
If $G$ is a group with a proper left-invariant metric,
then $\asdim(G)$ is the supremum of $\asdim(H)$ over all finitely generated subgroups $H$ of $G$.
\end{Cor}
\dokaz 
Since balls $B(1_G,r)$ are finite, the groups $G_r$ in  \ref{GroupAndSubgroups}
are finitely generated.
\hfill $\blacksquare$

\begin{Cor}[Bell-Dranishnikov \cite{BellDranish A Hurewicz Type}] \label{UnionThmForAsdim}
Suppose $X$ is a metric space and $\{X_s\}_{s\in S}$ is a family of subsets of $X$ such that
$X=\bigcup\limits_{s\in S}X_s$ and there is a single $n$-dimensional control function $D$ for all $X_s$. If, for every $r > 0$, there is a subspace $X_r$ of $X$
such that $\asdim(X_r)\leq n$ and the family $\{X_s\setminus X_r\}_{s\in S}$
is $r$-disjoint,
then $\asdim(X)\leq n$.
\end{Cor}
\dokaz 
Notice $\rdim(X\setminus X_r)\leq n$.
\hfill $\blacksquare$

\section{Dimension of a function}

When trying to generalize Hurewicz Theorem from covering dimension
to any other dimension theory, the issue arises of how to define
the dimension of a function. Let us present an example of
a Lipschitz function demonstrating that replacing
$\sup\{\dim(f^{-1}(y))\mid y\in Y\}$ by $\sup\{\asdim(f^{-1}(B))\mid B\subset Y \text{ is bounded}\}$
does not work. 

\begin{Prop} \label{ExampleOfWrongDimOfFunction}
There is a Lipschitz function $f\colon X\to Y$ of metric spaces 
such that $$\asdim(Y)=0=\sup\{\asdim(f^{-1}(B))\mid B\subset Y \text{ is bounded}\}$$
and $\asdim(X) > 0$.
\end{Prop}
\dokaz 
Let $Y$ consist of points $2^n$, $n\ge 1$, on the $x$-axis
and $X$ is the union of vertical segments $I_n$ of length $n$
and starting at $2^n$. $f\colon X\to Y$ is the projection.
Since $f^{-1}(B)$ is bounded for every bounded $B\subset Y$,
$\sup\{\asdim(f^{-1}(B))\mid B\subset Y \text{ is bounded}\}=0$.
Also, $\asdim(Y)=0$. However, $\asdim(X) > 0$
as for each $n$ it has $n$-components of arbitrarily large size.
\hfill $\blacksquare$

\begin{Rem}
\ref{ExampleOfWrongDimOfFunction} shows the answer to a problem
of Dranishnikov \cite{DranPrList} in negative. That problem asks
if $\Lasdim(X)\leq \Lasdim(Y)+\Lasdim f^{-1}$, where
$\Lasdim f^{-1}$ is defined as $\sup\{\Lasdim(f^{-1}(B))\mid B\subset Y \text{ is bounded}\}$
and $f$ is Lipschitz (see \ref{DefOfAssouadLinearlyControlled} for the definition
of $\Lasdim$).
Notice $\Lasdim(Y)=0$, $\Lasdim f^{-1}=0$, and $\Lasdim(X)\ge\asdim(X) > 0$
in \ref{ExampleOfWrongDimOfFunction}.
\end{Rem}

In view of \ref{ExampleOfWrongDimOfFunction} we define asymptotic
dimension of a function as follows:
\begin{Def} \label{DefOfAsDimOfFunction}
Given a function $f\colon X\to Y$ of metric spaces we define the {\it
asymptotic dimension} $\asdim (f)$ of $f$ as the supremum
of asymptotic dimensions of $A\subset X$ so that $f(A)\subset Y$ is
of asymptotic dimension $0$.
\end{Def}

\begin{Def} \label{DefOfDimControlOfFunction}
Given a function $f\colon X\to Y$ of metric spaces and given $m\ge 0$,
an  {\it $m$-dimensional control function} of $f$ is a function $D_f\colon R_+\times R_+\to R_+$
such that for all $r_X> 0$ and $R_Y > 0$ any $(\infty,R_Y)$-bounded subset $A$ of $X$ 
can be expressed as the
union of $m+1$ sets whose $r_X$-components are $D_f(r_X,R_Y)$-bounded.
\end{Def}

\begin{Prop} \label{DimOfFunctionAndUniformityOfDim}
Suppose $f\colon X\to Y$ is a function of metric spaces and $m\ge 0$.
If $\asdim (f)\leq m$, then $f$ has an $m$-dimensional control function $D_f$.
\end{Prop}
\dokaz Fix $r_X > 0$ and $R_Y > 0$.
Suppose for each $n$ there is $y_n\in Y$ such that $A_n=f^{-1}(B(y_n,R_Y))$
cannot be expressed as the
union of $m+1$ sets whose $r_X$-components are $n$-bounded.
The set $C=\bigcup\limits_{n=1}^\infty B(y_n,R_Y)$ cannot be bounded as $\asdim (f^{-1}(C))\leq m$
for any bounded subset $C$ of $Y$.
By passing to a subsequence we may arrange $y_n\to\infty$ and $\asdim (C)=0$, a contradiction.
\hfill $\blacksquare$

\begin{Def} \label{DefOfHigherDimControlOfFunction}
Given a function $f\colon X\to Y$ of metric spaces and given $k\ge m+1\ge 1$,
an  {\it $(m,k)$-dimensional control function} of $f$ is a function $D_f\colon R_+\times R_+\to R_+$
such that for all $r_X> 0$ and $R_Y > 0$ any $(\infty,R_Y)$-bounded subset $A$ of $X$ 
can be expressed as the
union of $k$ sets $\{A_i\}_{i=1}^k$ whose $r_X$-components are $D_f(r_X,R_Y)$-bounded
so that any $x\in A$ belongs to at least $k-m$ elements of $\{A_i\}_{i=1}^k$.
\end{Def}

\begin{Prop} \label{UniformityOfDimAndHigherCovers}
Let $f\colon X\to Y$ be a function of metric spaces and $m\ge 0$.
Suppose there is an $m$-dimensional control function $D^{(m+1)}_f\colon R_+\times R_+\to R_+$ of $f$. If one defines inductively functions $D^{(k)}_f$ for $k> m+1$ by $D^{(k)}_f(r_X,R_Y)=D^{(k-1)}_f(3r_X,R_Y)+2r_X$, then each $D^{(k)}_f$ is an $(m,k)$-dimensional
control function of $f$.
\end{Prop}
\dokaz The proof is  by induction on $k$. The
case of $k=m+1$ is obvious. Suppose the result holds for some $k \ge n+1$
and a subset $A$ of $X$ is $(\infty,R_Y)$-bounded.
There are $k$ subsets
$\{A_i\}_{i=1}^{k}$, each has $3r_X$-components bounded by
$D^{(k)}_f(3r_X,R_Y)$ such that the union of any $m+1$ of those sets covers $A$.

Define 
$A'_i$ to be the $r_X$-neighborhood of $A_i$ in $A$ for $i\leq k$.
Notice $r_X$-components of each $A'_i$ are $(D^{(k)}_f(3r_X,R_Y)+2r_X)$-bounded
as they are contained in $r_X$-neighborhoods of $3r_X$-components of $A_i$. 

Define
$A'_{k+1}$ as the union of all sets $\bigcap\limits_{i\in S}A_i\setminus \bigcup\limits_{i\notin S}A'_i$, where $S$ is a subset of $\{1,\ldots,k\}$ consisting of exactly
$k-m$ elements. 

Suppose $x\in A$ belongs exactly
to $k-m$ sets $A_i$ such that $i\leq k$ and let $S=\{i\leq k \mid x\in A_i\}$.
If $x\notin A'_{k+1}$, then $x$ must belong to $A'_j$ for some $j\notin S$.
Thus each $x\in A$ belongs to at least $k+1-m$ elements of $\{A'_i\}_{i=1}^{k+1}$.

Notice that any $r_X$-component of $A'_{k+1}$ is contained in a single
$r_X$-component of some $A_j$ resulting in
$r_X$-components of each $A'_i$ being $(D^{(k)}_f(3r_X,R_Y)+2r_X)$-bounded. 
\hfill $\blacksquare$

\begin{Prop} \label{UniformityOfDimAndDoubleComponents}
Let $f\colon X\to Y$ be a function of metric spaces and $k,m\ge 0$.
If $D_f\colon R_+\times R_+\to R_+$ is an $(m,k)$-dimensional control
function of $f$, then for any $B\subset Y$ whose $r_Y$-components
are $R_Y$-bounded, the set $f^{-1}(B)$
can be covered by $k$ sets whose $(r_X,r_Y)$-components are $D_f(r_X,R_Y)$-bounded
and every element of $f^{-1}(B)$ belongs to at least $k-m$ elements of that covering.
\end{Prop}
\dokaz Given an $r_Y$-component $S$ of $B$
express $f^{-1}(S)$ as $A_1^S\cup\ldots \cup A^S_k$
such that $r_X$-components of $A^S_i$ are $D_f(r_X,R_Y)$-bounded
and every element of $f^{-1}(S)$ belongs to at least $k-m$ elements of that covering.
Put $A_i=\bigcup\limits_S A^S_i$ and notice each $(r_X,r_Y)$-component
of $A_i$ is contained in an $r_X$-component of some $A^S_i$.
\hfill $\blacksquare$

\begin{Thm} \label{GeneralHurewicz}
Let $k=m+n+1$, where $m,n\ge 0$.
Suppose $f\colon X\to Y$ is a large scale uniform function of metric spaces and $\asdim (Y)\leq n$. 
$$\asdim (X)\leq m+n$$
 if there is an $(m,k)$-dimensional control function $D_f$ of $f$.
 Moreover, if one can choose the coarseness control function $c_f$ of $f$, the $n$-dimensional control function of $Y$, and $D_f$ to be linear (respectively, a dilation), then $X$ has a $(k-1)$-dimensional
 control function that is linear (respectively, a dilation).
\end{Thm}
\begin{pf} Let $c_f$ be a coarseness control function of $f$. Let $D_Y\colon R_+\to R_+$ be an $(n,k)$-dimensional control function
of $Y$. Notice we may require $c_f(r) > r$, $D_Y(r) > r$, and $D_f(r,R) > r+R$
as we may redefine those functions by adding $r$ or $r+R$ without losing their properties or type
(linear, dilation or dilation near $0$).

 Given a
number $r > 0$ we are going to construct a number $D_X(r)$ and represent
the space $X$ as a union of $k$ sets $\{D^j\}_{j=1}^k$ with all
$r$-components of $D^j$ being $D_X(r)$-bounded for every $j$.

Define inductively a sequence of numbers
$r^{(n+1)}_Y<R^{(n+1)}_Y<r^{(n)}_Y<R^{(n)}_Y<\dots<r^{(1)}_Y<R^{(1)}_Y<r^{(0)}_Y<R^{(0)}_Y$
starting from $r^{(n+1)}_Y=c_f(r)$ and moving to lower indices so that for every $i$ we have
$R^{(i)}_Y=D_Y(r^{(i)}_Y)$ and $r^{(i)}_Y=3R_Y^{(i+1)}$. 

Express $Y$ as the union of $k$ sets $\{A_i\}_{i=1}^{k}$
such that all $r^{(0)}_Y$-components of $A_i$ are
$R^{(0)}_Y$-bounded for every $i$. Since $A_i\subset Y$, we can express the set $A_i$ as the union of $k$ sets
$\{U_i^j\}_{j=1}^{k}$ such that all $r^{(i)}_Y$-components of
$U_i^j$ are $R^{(i)}_Y$-bounded for every $j$ and every point
$y\in A_i$ belongs to at least $m$ sets.

Define inductively a sequence of numbers
$r^{(n+1)}_X<R^{(n+1)}_X<r^{(n)}_X<R^{(n)}_X\dots<r^{(1)}_X<R^{(1)}_X$
starting with
 $r^{(n+1)}_X=r$ and for every $i$
we have $R^{(i)}_X=D_f(r^{(i)}_X,R^{(i)}_Y)$ and
$r^{(i)}_X=3R_X^{(i+1)}$. 

For every $i$ we express the set $f^{-1}(A_i)$ as the union of $k$ sets
$\{B_i^j\}_{j=1}^{k}$ such that all $(r^{(i)}_X,r^{(i)}_Y)$-components of
$B_i^j$ are $(R^{(i)}_X,R^{(i)}_Y)$-bounded for every $j$ and every point
$x\in f^{-1}(A_i)$ belongs to at least $n$ sets.

Put $D^j_i=B^j_i\cap f^{-1}(U^j_i)$ and let $D^j$ be the union of
all $D^j_i$. Notice $D^j$'s cover $X$ by the use of Kolmogorov's
argument: given $x\in X$ there is $i$ so that $f(x)\in A_i$. The
set of $j$'s such that $x\in B^j_i$ has at least $k-m$ elements,
the set of $j$'s such that $f(x)\in U^j_i$ has at least $k-n$
elements, so they cannot be disjoint.

Notice all
$(r^{(i)}_X,r^{(i)}_Y)$-components of the set $D_i^j$ are
$(R^{(i)}_X,R^{(i)}_Y)$-bounded. By
\ref{rComponentsOfBigUnions} all
$(r^{(n+1)}_X,r^{(n+1)}_Y)$-components of the set $D^j$ are
$(3R^{(1)}_X,3R^{(1)}_Y)$-bounded. 
Since $(r^{(n+1)}_X,r^{(n+1)}_Y)=(r,c_f(r))$,
by \ref{CoarsnessLipViaComponents} all
$r$-components of
$D^j$ are $3R^{(1)}_X$-bounded.
Observe, $3R^{(1)}_X$ is a linear function of $r$ (respectively, a dilation), if
the coarseness control function $c_f$ of $f$, the $n$-dimensional control function of $Y$, and $D_f$ are linear (respectively, dilations).
\end{pf}

\begin{Cor} \label{DimOfFunctionAndUniformityOfDimEquiv}
Suppose $f\colon X\to Y$ is a function of metric spaces and $m\ge 0$.
$\asdim (f)\leq m$ if and only if $f$ has an $m$-dimensional control function.
\end{Cor}
\dokaz In one direction use \ref{DimOfFunctionAndUniformityOfDim}.
In the other direction apply \ref{GeneralHurewicz} in the case of $n=0$.
\hfill $\blacksquare$

\begin{Thm} \label{HurForAsdim}
$$\asdim (X)\leq \asdim (f)+\asdim (Y)$$ for any large scale uniform
function $f\colon X\to Y$.
\end{Thm}
\dokaz Follows from \ref{GeneralHurewicz}, \ref{DimOfFunctionAndUniformityOfDimEquiv}, and \ref{UniformityOfDimAndHigherCovers}.
\hfill $\blacksquare$

In \cite{BellDranish A Hurewicz
Type} there is a concept of a family $\{X_{\alpha}\}$
of subsets of $X$ satisfying $\asdim(X_{\alpha}) \le n$ uniformly. Notice that in our language
it means there is one function that serves as an $n$-dimensional control function
for all $X_{\alpha}$.

\begin{Cor}[Bell-Dranishnikov \cite{BellDranish A Hurewicz
Type}]\label{HurewiczTypeOfBellD}
Let $f\colon X \to Y$ be a Lipschitz function of metric spaces. Suppose that,
for every $R >0$, 
$$\asdim\{f^{-1}(B_R (y))\} \le n$$
 uniformly
(in $y \in Y$). If $X$ is geodesic, then $\asdim(X) \le \asdim(Y) +n$.
\end{Cor}
\dokaz
$\asdim\{f^{-1}(B_R (y))\} \le n$ uniformly means $f$ has an $n$-dimensional
control function, so apply  \ref{GeneralHurewicz} and \ref{UniformityOfDimAndHigherCovers}.
\hfill $\blacksquare$

\begin{Rem} Notice
Bell-Dranishnikov's version of Hurewicz type theorem \ref{HurewiczTypeOfBellD}
(see \cite{BellDranish A Hurewicz Type}) assumes that $X$ is  geodesic. We
do not use this assumption in \ref{HurForAsdim}.
\end{Rem}

\section{Asymptotic dimension of groups}
J.Smith \cite{Smith} showed that any two proper and left-invariant metrics on
a given countable group $G$ are coarsely equivalent.
We generalize that result as follows.

\begin{Prop} \label{HomomorphismsAreCoarseIfCoarselyProper}
Suppose $f\colon G\to H$ is a homomorphism
of groups and $d_G$, $d_H$ are left-invariant metrics on $G$ and $H$, respectively.
If $f\colon (G,d_G)\to (H,d_H)$ is coarsely proper (i.e., it sends bounded subsets of $G$ to bounded subsets of $H$), then $f\colon (G,d_G)\to (H,d_H)$ is large scale uniform.
\end{Prop}
\dokaz 
Suppose $r > 0$. Since $f(B(1_G,r))$ is bounded, there is $R > 0$ such that $d_H(1_H,f(g)) < R$ for all $g\in B(1_G,r)$. If $x,y\in G$ satisfy $d_G(x,y) < r$, then $x^{-1}\cdot y\in B(1_G,r)$,
so $d_H(f(x),f(y))=d_H(1_H,f(x^{-1}y)) < R$.
\hfill $\blacksquare$

\begin{Cor} \label{HomomorphismsAreCoarse}
Suppose $f\colon G\to H$ is a homomorphism
of groups and $d_G$, $d_H$ are left-invariant metrics on $G$ and $H$, respectively.
If $d_G$ is proper, then $f\colon (G,d_G)\to (H,d_H)$ is large scale uniform.
\end{Cor}
\dokaz 
Since bounded subsets of $G$ are finite, $f$ is coarsely proper.
\hfill $\blacksquare$

\begin{Cor}[J.Smith \cite{Smith}] \label{HomomorphismsAreCoarse}
Any two proper left-invariant metrics on a group $G$ are coarsely equivalent.
\end{Cor}

The subsequent result is derived in \cite{DS} from the corresponding theorem for finitely generated groups. We deduce it directly from \ref{HurForAsdim}.
\begin{Thm}[Dranishnikov-Smith \cite{DS}]\label{AsDimOfExactSequence}
If $1\to K\to G\to H\to 1$ is a short exact sequence of countable groups, then
$$\asdim(G)\leq \asdim(K)+\asdim(H).$$
\end{Thm}
\dokaz 
It suffices to show $\asdim(f)=k=\asdim(K)$, where $f\colon G\to H$. 
Given $R_H > 0$, the subspace $f^{-1}(B(1_H,R_H))$ is in a bounded neighborhood
of $K$ as $B(1_H,R_H)$ is finite, so it has the same asymptotic dimension 
as $K$. Fix the dimension control function $D(r_G,R_H)$ for that space.
Notice that all spaces $f^{-1}(B(h,R_H))$, $h\in H$, are isometric to $f^{-1}(B(1_H,R_H))$,
so $D(r_G,R_H)$ is a dimension control function of all of them.
By \ref{DimOfFunctionAndUniformityOfDimEquiv}, $\asdim(f)\leq k$ and by applying 
\ref{HurForAsdim} we are done.
\hfill $\blacksquare$

In \cite{DS} the asymptotic dimension of a group $G$ is defined as the supremum of
$\asdim(H)$, $H$ ranging through all finitely generated subgroups of $G$.
We provide an alternative definition which will be applied in the next section.

\begin{Def} \label{DefOfSComponentsInAnyGroup}
Given a finite subset $S$ of a group $G$, an {\it $S$-component} of $G$
is an equivalence class of $G$ of the relation $x\sim y$ iff
$x$ can be connected to $y$ by a finite chain $x_i$ so that $x_i^{-1}\cdot x_{i+1}\in S$
for all $i$.

A family of subsets $\{A_i\}_{i\in J}$ is {\it $S$-bounded} if $x^{-1}\cdot y\in S$ for each $i\in J$ and all $x,y\in A_i$.
\end{Def}

\begin{Prop}\label{AsDimViaFiniteSets}
Let $G$ be a group. $\asdim(G)\leq n$ if and only if for each finite subset $S$ of $G$
there is a finite subset $T$ of $G$ and a decomposition $G=A_1\cup\ldots\cup A_{n+1}$
such that $S$-components of $A_i$ are $T$-bounded for all $i\leq n+1$.
\end{Prop}
\dokaz 
Suppose for each finite subset $S$ of $G$
there is a finite subset $T$ of $G$ and a decomposition $G=A_1\cup\ldots\cup A_{n+1}$
such that $S$-components of $A_i$ are $T$-bounded for all $i\leq n+1$.
Given a countable subgroup $H$ of $G$ and given a proper left-invariant
metric $d_H$ on $H$, we put $S=B(1_H,r)$ and $R=\max\{d_H(1_H,t) \mid t\in T\cdot T\cap H\}$ (by $T\cdot T$ we mean all products $t\cdot s$, where $t,s\in T$).
Notice that $r$-components of $A_i\cap H$ are $R$-bounded.
\par Conversely, suppose $\asdim(H)\leq n$ for all finitely generated subgroups $H$ of $G$.
Given a finite subset $S$ of $G$ let $H$ be the subgroup of $G$ generated by $S$
and let $d$ be the word metric on $H$ induced by $S$. Since $\asdim(H)\leq n$,
there is a decomposition of $H$ into $A_1\cup\ldots \cup A_{n+1}$ such that
$1$-components of $A_i$ are $m$-bounded for all $i\leq n+1$ and some $m > 0$.
Let $T=B(1_H,m+1)$. Pick a representative $g_j\in G$, $j\in J$, of each left coset
of $H$ in $G$. Put $B_i=\bigcup\limits_{j\in J}g_j\cdot A_i$ and notice $S$ components of $B_i$
are $T$-bounded.
\hfill $\blacksquare$

\section{Linearly controlled asymptotic dimension of groups}

Unlike the asymptotic dimension of countable groups, the
asymptotic  Assouad-Nagata dimension may depend on the left-invariant metric.
Piotr Nowak \cite{Nowak} constructed finitely generated groups $G_n$
of asymptotic dimension $n\ge 2$ and $\ANasdim (G_n)=\infty$.
Obviously, there is no such example for $n=0$ as finitely generated groups
of asymptotic dimension $0$ are finite. However, \ref{MetricOnTorsionWithDimInfty} does provide
a countable group $G$ with a proper, left-invariant metric $d$ such that
$\asdim(G,d)=0$ and $\ANasdim(G,d)=\infty$.

The purpose of this section is to solve in negative
the following two problems of Dranishnikov \cite{DranPrList}.
\begin{Problem}\label{Problem 37}
Does $\Lasdim(X)=\ANasdim(X)$ hold for metric spaces?
\end{Problem}

\begin{Problem} \label{Problem 41c}
Find a metric space $X$ of minimal
$\asdim(X)$ such that $\asdim(X) < \Lasdim(X)$.
\end{Problem}

Here is an answer to  \ref{Problem 41c}.

\begin{Prop} \label{MetricOnTorsionWithLDimBiggerThan0}
There is a proper, left-invariant metric $d_G$ on $G=\bigoplus\limits_{n=2}^\infty \ZZ /n$
such that $\Lasdim(G,d_G) > 0=\asdim(G,d_G)$.
\end{Prop}
\dokaz Pick a generator $g_n$ of $\ZZ /n$ and assign it the norm of $n$
for $n\ge 2$.
Extend the norm over all elements $g$ of $G$ as the minimum of $\sum\limits_{n=1}^\infty |k_n|\cdot n$,
where $g=\sum\limits_{n=1}^\infty k_n\cdot g_n$ and $k_n$ are integers.
Notice that $d_G(g,h):=|g-h|$ is a proper and invariant metric on $G$.
Suppose there is an unbounded subset $U$ of $\RR_+$
such that for some $C > 0$ all $r$-components of $G$ are
$(C\cdot r)$-bounded for $r\in U$. Pick $r\in U$ such that $4C+4 < r$ and let $n$ be an integer
such that $r-2< 2n \leq r$. Notice all $2n$-components of $G$ are contained
in some $r$-components of $G$, so they are $(n-1)\cdot (n+1)$-bounded
as $C < \frac{r-4}{4} \leq \frac{n-2}{2}$ and $r < 2n+2$.
However, the $2n$-component of $0$ contains $\ZZ/2n$ and $|n\cdot g_{2n}|=2n^2 > (n-1)\cdot (n+1)$,
a contradiction.
\hfill $\blacksquare$

\begin{Prop} \label{MetricOnTorsionWithDimInfty}
There is an Abelian torsion group $G$ with a proper invariant metric $d_G$  such that 
$\Lasdim (G,d_G) =0$ and $\ANasdim (G,d_G) =\infty$.
\end{Prop}
\dokaz 
Consider $\ZZ^{n+1}$ with the standard word metric $\rho_n$. Since $\asdim(\ZZ^{n+1})> n$,
there is $r(n) > 0$ such that any decomposition of $\ZZ^{n+1}$ as
the union $X_0\cup\ldots\cup X_n$ forces one $X_i$ to have $r(n)$-components
without a common upper bound.
Without loss of generality we may assume $r(n)$ is an integer greater than $n$.
Put $s(n)=8\cdot n\cdot r(n)$ and consider $G_n=(\ZZ_{s(n)})^{n+1}$ with
the metric $d_n$ equal to $t(n)$ times the standard word metric.
Numbers $t(n)$ are chosen so that $t(n+1) > \diam(G_{1})+\ldots+\diam(G_{n})$
for all $n\ge 1$. Since $\diam(G_{n})=4nt(n)r(n)$,
we want $t(n+1) =1+ 4\cdot t(1)r(1)+\ldots+4n\cdot t(n)r(n)$ for all $n\ge 1$.

The group $G$ is the direct sum of all $G_n$ with the obvious metric $d$:
$d(\{x_n\},\{y_n\})=\sum d_n(x_n,y_n)$.
Notice that $\frac{t(n)}{2}$-components of $G$ are of size at most $t(n)$.
Indeed, any $\frac{t(n)}{2}$-component of $G$ 
is contained in $x+G_1\oplus\ldots\oplus G_{n-1}$ for some $x\in G$.
Thus $\Lasdim(G,d)=0$.
\par

Consider the projection $\pi_n\colon \ZZ^{n+1}\to G_n$ and notice it is $t(n)$-Lipschitz.
Moreover, if two points $x$ and $y$ in $\ZZ^{n+1}$ satisfy
$\rho_n(x,y)\leq 4nr(n)$, then $d_n(x,y)=t(n)\cdot \rho_n(x,y)$.

Suppose $G_n$ equals $Y_0\cup\ldots \cup Y_n$. Assume $r(n)$-components
of some $X_j=\pi_n^{-1}(Y_j)$ are not uniformly bounded. Therefore there is a sequence
$x_0,\ldots,x_k$ in $X_j$ such that $\rho_n(x_i,x_{i+1})< r(n)$
and $(n+2)\cdot r(n) > \rho_n(x_0,x_k)> n\cdot r(n)$.
Assume such a sequence does not exist with number of elements smaller than $k+1$.
Therefore $d_n(\pi_n(x_i),\pi_n(x_{i+1}))=t(n)\cdot \rho_n(x_i,x_{i+1}) < t(n)\cdot r(n)$
and $(n+2)\cdot t(n)r(n) > \rho_n(\pi_n(x_0),\pi_n(x_k))> n\cdot t(n)r(n)$.
Thus, one $(t(n)\cdot r(n))$-component of $Y_i$ is of diameter bigger than $n\cdot t(n)r(n)$.
\par Suppose $\ANasdim (G,d)=k < \infty$ and $D_G(r)=c\cdot r+b$ is a $k$-dimension
control function of $(G,d)$. Pick $n > \max(k,|c|,|b|)+2$ and choose a decomposition
$Y_0\cup\ldots\cup Y_k$ of $G$ so that $(t(n)\cdot r(n))$-components of each $Y_i$
are bounded by $c\cdot t(n)\cdot r(n)+b < (c\cdot t(n)+1)\cdot r(n)$. Notice there is a $1$-Lipschitz
projection $p_n\colon G\to G_n$. Therefore each set $p_n(Y_i)$ has $(t(n)\cdot r(n))$-components
bounded by $(c\cdot t(n)+1)\cdot r(n)$. On the other hand,
for some $j$, there is an $(t(n)\cdot r(n))$-component of $p_n(Y_j)$ of diameter bigger than
$nt(n)\cdot r(n)$. Thus, $c\cdot t(n)+1> nt(n)$ and $c > n-1/t(n)$, a contradiction.
\hfill $\blacksquare$

\begin{Rem} Proposition \ref{MetricOnTorsionWithDimInfty} solves Problem
\ref{Problem 37} in negative.
\end{Rem}

\begin{Prop} \label{LasdimMayEqualasdim}
For any countable group $G$
there is a proper left-invariant metric $d_G$ such that $\ANasdim (G,d_G)=\asdim(G)$.
\end{Prop}
\dokaz Let $\asdim(G)=n-1$.
Express $G$ as the union of an increasing sequence $S_0\subset S_1\ldots$
of its finite subsets so that $\{1_G\}=S_0$, each $S_i$ is symmetric,
and for each $i$ there is a decomposition of $G$ as $A^i_1\cup\ldots \cup A^i_n$
so that $S_i$-components of $A^i_k$ are $S_{i+1}$-bounded.
Define $d_G(x,y)$ as the smallest $i$ so that $x^{-1}\cdot y\in S_i$.
Notice that in the metric $d_G$ all $i$-components of $A^i_k$ are $(i+1)$-bounded.
\hfill $\blacksquare$

\section{Assouad-Nagata dimension}

\begin{Def} \label{DefOfNADimOfFunction}
Given a function $f\colon X\to Y$ of metric spaces we define the {\it
Assouad-Nagata dimension} $\dim_{AN}(f)$ of $f$ as the minimum
of $m$ for which there is an $m$-dimensional control function $D_f(r_X,R_Y)$
of the form $a\cdot r_X+b\cdot R_Y$.
\end{Def}

\begin{Thm} \label{HurThmForDilationControlled}
If $f\colon X\to Y$ is a Lipschitz function of metric spaces, then
$\dim_{AN}(X)\leq \dim_{AN} (f)+\dim_{AN}(Y)$.
\end{Thm}
\dokaz 
Notice in the proof of \ref{GeneralHurewicz} the resulting dimensional control function of $X$
is a dilation if $c_f$, $D_X$, and $D_f$ are dilations.
\hfill $\blacksquare$

\begin{Rem}  Notice the proof of of \ref{GeneralHurewicz} can be analyzed to give a version of \ref{HurThmForDilationControlled} for microscopic Assouad-Nagata dimension.
\end{Rem}

\section{Asymptotic Assouad-Nagata dimension}

\begin{Def} \label{DefOfLinAsDimOfFunction}
Given a function $f\colon X\to Y$ of metric spaces we define the {\it
asymptotic  Assouad-Nagata dimension} $\ANasdim  (f)$ of $f$ as the minimum
of $m$ for which there is an $m$-dimensional control function $D_f(r_X,R_Y)$
of the form $a\cdot r_X+b\cdot R_Y+c$.
\end{Def}

\begin{Thm} \label{HurThmForLinearlyControlled}
If $f\colon X\to Y$ is an asymptotically Lipschitz function of metric spaces, then
$\ANasdim (X)\leq \ANasdim  (f)+\ANasdim (Y)$.
\end{Thm}
\dokaz 
Notice in the proof of \ref{GeneralHurewicz} the resulting dimensional control function of $X$
is linear if $c_f$, $D_X$, and $D_f$ are linear.
\hfill $\blacksquare$

\begin{Prop} \label{LAsdimForCountableGroups}
Let $n\ge 0$. If $(G,d_G)$ is a group equipped with a proper, left-invariant metric $d_G$,
then the following conditions are equivalent:
\begin{itemize}
\item[a.] $\ANasdim (G,d_G)\leq n$.
\item[b.] There are constants $C, M > 0$ such that the function $r\to M\cdot r+C$
is an $n$-dimensional control function for all finitely generated subgroups
of $G$.
\end{itemize}
\end{Prop}
\dokaz a)$\implies$b) follows from the proof of \ref{GroupAndSubgroups}.
b)$\implies$a) is obvious.
\hfill $\blacksquare$

\begin{Prop} \label{MetricsForEpi}
If $1\to K\to G\to H\to 1$ is an exact sequence and $G$ is a finitely generated group, then
there are word metrics $d_G$ on $G$ and $D_H$ on $H$
such that $f\colon (G,d_G)\to (H,d_H)$ is $1$-Lipschitz and for any $m$-dimensional control function $D_K$ on $K$
the function $$D_f(r_G,R_H):=D_K(r_G+2R_H)+2R_H$$
is an $m$-dimensional control function of $f$.
\end{Prop}
\dokaz 
Let $S$ be a symmetric finite set of generators for $G$.
Let $d_G$ be the word metric on $G$ induced by $S$ and let
$d_H$ be the word metric on $H$ induced by $f(S)$.
\par
Let $B=B(1_H,R_H)$ and $A=f^{-1}(B)$. Notice $A\subset B(K,R_H)$.
Indeed, if $d_H(1_H,f(a)) < i=R_H$, then $f(a)=\prod\limits_{i=1}^mf(s_i)$
such that $k < R_H$ and $s_i\in S$ for $i\leq m$. Consequently, $a=k\cdot x$ so that $k\in K$
and $x=\prod\limits_{i=1}^ms_i$ has the property $d_G(1_G,x)< R_H$.

Now $d_G(k\cdot s,k)=d_G(s,1_G)\leq i-1< i$ and $a\in B(K,R_H)$.
By \ref{DimControlOfAneighborhood}, the function $D_f(r_G,R_H):=D_K(r_G+2R_H)+2R_H$
is an $m$-dimensional control function of $B(K,R_H)$.
Since all $f^{-1}(B(y,R_H))$ are isometric to $A$, we are done.
\hfill $\blacksquare$

\begin{Cor} \label{HurThmForLinearlyControlledExactSeqOfFinGenerated}
If $1\to K\to G\to H\to 1$ is an exact sequence of groups so that $G$ is finitely generated, then
$$\ANasdim (G,d_G)\leq \ANasdim  (K,d_G|K)+\ANasdim (H,d_H)$$
for any word metrics metrics $d_G$ on $G$ and $d_H$ on $H$.
\end{Cor}
\dokaz 
Use the metrics as in \ref{MetricsForEpi}. That way 
$$\ANasdim (f)\leq \ANasdim  (K,d_G|K),$$
so applying \ref{HurThmForLinearlyControlled} one gets the desired inequality.
\hfill $\blacksquare$

One is tempted to define the linear asymptotic dimension of arbitrary groups
as the supremum of $\ANasdim (H)$ for all finitely generated subgroups $H$ of $G$,
However, one runs into problems with that definition.

\begin{Question}\label{IsLasdimmonotonic}
Suppose $G$ is a finitely generated group and $H$ is its finitely generated subgroup.
Does $\ANasdim (H)\leq \ANasdim (G)$ hold? Is there a case of $\ANasdim (H)=\infty$
and $\ANasdim (G) < \infty$? 
\end{Question}

\begin{Question}\label{IsLasdimAsSubgroupEqualLasdim}
Suppose $G$ is a finitely generated group and $H$ is its finitely generated subgroup.
Does $\ANasdim (H)=\ANasdim (H,d_G|H)$ hold for any word metric $d_G$ on $G$?
\end{Question}

\begin{Question}
Suppose $G$ and $H$ are two finitely generated groups of finite
asymptotic  Assouad-Nagata dimension. Is $\ANasdim (G\ast H)$ finite?
Does $\ANasdim (G\ast H)=\max(\ANasdim (G),\ANasdim (H),1)$ hold?
\end{Question}

To complete this section we state and prove a version of the the Hurewicz theorem for
groups acting on spaces of finite asymptotic  Assouad-Nagata dimension. For
this we need the concept of $R$-stabilizers. Let $G$ act on the
metric space $X$ by isometries and let $R > 0$. Given $x_0 \in X$
the {\it $R$-stabilizer} of $x_0$ is defined by $W_R(x_0) = \{\gamma \in
G: d(\gamma\cdot x_0,x_0) \le R\}$.

\begin{Thm}\label{HurewiczforFGGroupsactingonSpcs}
Let $G$ be a finitely generated group acting by isometries on a metric
space $X$ of finite asymptotic  Assouad-Nagata dimension. Fix
a point $x_0 \in X$. If there are constants $m,b,c > 0$ such that
$r\to m\cdot r+b\cdot R+c$ is a $k$-dimensional control function of
$W_R(x_0)$ for
all $R$, then $$\ANasdim (G) \le k + \ANasdim (X).$$
\end{Thm}

\dokaz Let $S$ be a symmetric generating set for $G$ and let
$\lambda=\max\{d_X(s\cdot x_0,x_0)|s \in S\}$. Define $\pi\colon G \to X$
by $\pi(\gamma)=\gamma\cdot x_0$ and notice that $\pi$ is
$\lambda$-Lipschitz. Also,
${\pi}^{-1}(B_R(g\cdot x_0))=gW_R(x_0)$ shows that the sets
$({\pi}^{-1}(B_R(g\cdot x_0)))$ are all isometric to $W_R(x_0)$ and we can use
\ref{HurThmForLinearlyControlled} to get our result.
\hfill $\blacksquare$

\end{document}